\documentclass[twoside,reqno,12pt]{article} %%{fcaa-var} %

 \usepackage{hyperref} % Editor will use to create hyperlinks %
\usepackage{amsmath}

\usepackage{amssymb}
\usepackage{graphicx}

\usepackage{mathrsfs}

\usepackage{color}

\newcommand{\be}{\begin{eqnarray}}
\newcommand{\ee}{\end{eqnarray}}

\def\eg{{\it e.g. }} 
\def\ie{{\it i.e. }}

\def\RR{\vbox {\hbox to 8.9pt {I\hskip-2.1pt R\hfil}}}

\def\CC{{\rm C\hskip-4.8pt \vrule height 6pt width 12000sp\hskip 5pt}}

\def\ds{\displaystyle}
\def\q{\quad} 
\def \rec#1{\frac{1}{#1}}
\def\ds{\displaystyle}
\def\eg{{\it e.g.}\ }
\def\ie{{\it i.e.}\ }

\def\e{{\rm e}}
 
 \def\L{{\mathcal{L}}}

\def\e{{\rm e}}

\def\Ei{{\rm Ei}\,}
\def\Ein{{\rm Ein}\,}
\def\log{{\rm log}\,}
\def\EE{{\mathcal E}}

\font\title=cmbx12 scaled\magstep2

\font\little=cmr10
\begin{document}
\begin{center}
%%{\bfs Time-fractional derivatives in relaxation processes:}
%%\vs{\bfs a tutorial survey}\vvs
%%%%%%%%%%%%%%%%%%%%%%%%%%%%%%%%%%%%%%%%%%%%%%%%%%%%%%%%%%%%%%%%%%%%%%%%%
{\title {On modifications of the exponential}}  
\\ [0.25truecm]
{\title{ integral with the Mittag-Leffler function}}
%% \footnote{This paper would be referred to as
%% {\it Eur. Phys. J. Special Topics} {\bf 193}, 133--160 (2011).}}
\\ [0.25truecm]
{Francesco MAINARDI} $^{(1)}$ and
{Enrico MASINA}$^{(2)}$
%%%%%%%%%%%%%%%%
\\
$\null^{(1)}$
 {\little Department of Physics and Astronomy, University of Bologna, and INFN,} \\
              %%%  Sezione di Bologna, \\
{\little Via Irnerio 46, I-40126 Bologna, Italy} \\
%% Tel: +39-051.2091098 $\;$ Fax: +39-051.247244 $\;$\\
{\little Corresponding Author. E-mail: francesco.mainardi@ibo.infn.it} 
\\ [0.25 truecm]
$\null^{(2)}$
{\little  Department of Physics and Astronomy, University of Bologna, and INFN,} \\
{\little  Via Irnerio 46, I-40126 Bologna, Italy}
\end{center}

% \vbox to 2.5cm { \vfill }

%%% to make empty space of approx. 2.5cm %%%%%%
%%% will be replaced by Editor with the journal's and publoishers logos %%%%%%%%

 %\bigskip \medskip

%%%% Abstract %%%%%%%%%%%%%%%%%%%%%%%%% 

\begin{abstract}
In this paper we survey the properties of the  
Schelkunoff modification of the Exponential integral and we generalize it with the Mittag-Leffler function.
So doing we get a new special function (as far as we know)  that may be  relevant in linear viscoelasticity because of its complete monotonicity properties in the time domain.   
We also consider the generalized sine and cosine integral functions.
 \end{abstract}
 
%%%%


\smallskip
\noindent
{\bf Paper published in
Fract. Calc. Appl. Anal (FCAA)., Vol. 21, No 5 (2018), pp. 1156--1169.
{DOI: 10.1515/fca-2018-0063};\\
ERRATUM in FCAA, Vol 23, No 2 (2020), pp 600--603. 
{DOI: 10.1515/fca-2020-0030}.
}
 
 %}\end{abstract}

 %\maketitle

%%%%%%%%%%%%%%%
\section*{Foreword for the Revised Version}
This is a revised version taking into account the  correct plot in the right of Fig 4
in the previous version.
It has been observed by Richard Paris \cite{Paris_2020} that the plot concerning the generalized cosine integral in Fig 4-Right is incorrect.
As a consequence we now provide the correct plot. By the way we take this occasion to pint out the article by Paris \cite{Paris_2020} where the author has provided a detailed study of the asymptotic expansion of the functions introduced in his paper that reduce as particular cases those dealt herewith. 
 
\section{Introduction}
One purpose of this paper is to point out the relevance of a modification due to  Schelkunoff on 1944
of the Exponential integral that turns out to 
be not sufficiently pointed out
%5 be practically ignored 
in the majority of textbooks on special functions.
The merit of this modification was later recognized by the great Italian mathematician F.G Tricomi
because it is proved to  provide an entire function. Here we propose to generalize this modification by introducing in the kernel of the integral the Mittag-Leffler function which is known to be the most direct generalization of the exponential function like it is the gamma function for the factorial.

The plane of the paper is the following.
In Section 2 we recall the known expressions of the Exponential integral and introduce the Schelkunoff modification. This modification leads to an entire  function in the complex plane 
Ein(z)
that in the time positive domain $t>0$ is of Bernstein type, namely its derivative is completely monotone. This property is equivalent to have a non-negative spectral distribution    that indeed was adopted by Becker in 1926 for an interesting model in linear viscoelasticity even if independently of the later  note by Schelkunoff.

In Section 3 we introduce the Mittag-Leffler function depending on a positive parameter $\nu$ in the kernel in such a way that for $\nu=1$ the transcendental  function $Ein_\nu(z)$  reduces to the Schelkunoff function and for $0<\nu\le 1$ and  $t>0$ keeps the Bernstein character thanks the well-known complete monotonicty of the Mittag-Leffler function. We consider also the limiting case $\nu=0$ and  recall for $0\le\nu\le 1$
the relevance of our function in describing the creep features of a viscoelastic model already discussed recently by Mainardi et al \cite{Mainardi-Masina-Spada_2019}.
For our function we exhibit the spectral distributions for selected values of $\nu \in (0,1)$  derived from  its Laplace transform by using the Mathematica tool box.  

In section 4, after recalling the  Schelkunoff modification of the functions related to the exponential integrals: sine and  cosine exponential integrals,  we generalize them by introducing in the kernel of the integral the Mittag-Leffler function. 
On this respect, we  adopt for fractional circular functions the definitions by 
Herrmann \cite{Herrmann_BOOK2014} noting that there exist  
different definitions for them in the literature.
% two alternative and different approaches in absence of a unique definition of fractional circular functions accepted in the literature. 
% Then we show for consistency our definitions of fractional sine and cosine % functions. 
  
Finally, in Section 5 we draw our conclusions. 

\section{The Exponential integral and its Schelkunoff modification}
  %$ \Ei(z)$, $\EE_1(z)$}

\subsection{The classical Exponential integral $\Ei(z)$}%%
A classical definition of the {\it Exponential integral} is
%% $$\E_n(z) :=
%%  \int_1^\infty \!\! \frac{\;\e^{\ds -zt}}{{\ds t^n}}\, dt $$
$$  \Ei(z) :=
 - \int_{-z}^\infty \!\! \frac{\;\e^{\ds-u}}{{\ds u}}\, du \,, 
  \q z\in \CC^-\,,
    \eqno(2.1)
$$
where the Cauchy principal  value of the integral is understood
if $z=x>0$.
%%%%%%
Above and from now on,
with $\CC^-$ we denote the complex plane $\CC$ cut along the 
negative real negative axis, that is $|\hbox{arg}\, z|<\pi$.

Some authors such as Jahnke and Emde \cite{Jahnke-Emde_BOOK1943} adopt the following definition
for $\Ei(z)\,, $
 $$\Ei(z) :=  \int_{-\infty}^{z} \!\!
 \frac{\;\e^{u}}{{\ds u}}\, du\,,
  \q z\in \CC^-. 
    \eqno (2.2)$$
which is equivalent to (2.1). For this we note that
 $$ \int_{-z}^\infty \!\! \frac{\;\e^{-u}}{{\ds u}}\, du
       = -   \int_{-\infty}^{z} \!\!
 \frac{\;\e^{u}}{{\ds u}}\, du
\,,
 $$
where the Cauchy principal value is understood for $z=x>0\,. $

Recalling  the  incomplete Gamma functions  
$$\Gamma(\nu, z) := \int_z^\infty u^{\nu-1} \e^{\ds -u}\, du\,,  \q \nu \in \RR\,, \q 
 z\in \CC^-, $$
we note
the  {identity}
$$   -\Ei(-z) = \Gamma(0\,,z) =
     \int_{z}^\infty \!\! \frac{\;\e^{\ds-u}}{ {\ds u}}\, du
   \,, \q z\in \CC^-\,.    \eqno(2.3)     $$

\subsection{The function $\EE(z)$}
%%%
In many texts on special functions the function $\Gamma(0\,,z) $
is usually taken as definition of Exponential integral and denoted  by
$\EE_1(z)$.  Then  we can  account for the following equivalent expressions
for  $z \in \CC^-$: 
$$ \EE_1(z) =  \Gamma(0\,,z)=  -\Ei(-z)  =
     \int_{z}^\infty \!\! \frac{\;\e^{\ds-u}}{ {\ds u}}\, du
   = \int_1^\infty \frac{\e^{-zt}}{t} \, dt \,.\eqno(2.4)$$
This definition is then generalized to yield, see \eg \cite{NIST}
$$ \EE_\nu(z)=  \int_1^\infty \frac{\e^{-zt}}{t^\nu } \, dt  =
z^{\nu-1}\, \Gamma (1-\nu, z) \,,\q
 \nu \in \RR\,,
\q z \in \CC^- 
  \,. \eqno(2.5)$$
%%%%%%%%%
We note that, in contrast with  the standard literature
where the Exponential integrals
are  denoted by the letter $E$,
we have used for them the letter $\EE$:
this  choice is to   avoid  confusion with the
standard notation for the Mittag-Leffler function
$E_\nu (z)$ ($\nu>0$)  later used in this paper to generalize the Exponential integral.
  We recall that the Mittag-Leffler function  plays a relevant role in fractional calculus; 
for more details see any treatise on fractional calculus and in particular the 2014 monograph by Gorenflo et al.
\cite{GKMR_BOOK2014}.
%%%%%%%%%%%%%%%%%%%

\subsection{The modified Exponential integral $\Ein(z)$}
%
%\subsection{The basic definition.}
The whole subject matter can be greatly simplified if we agree to
follow F.G. Tricomi \cite{Tricomi_BOOK1959}, see also his former assistant L. Gatteschi \cite{Gatteschi_BOOK1973},    who has proposed to consider
the following {\it entire}  function, formerly introduced in 1944 by Schelkunoff \cite{Schelkunoff_1944}:
%% ["Proposed symbols for the modified cosine and integral exponential integral",
%% {\it Quart. Appl. Math.} {\bf 2} (1944), p. 90]:
$$ \Ein(z) :=  \int_0^z \frac{1-\e^{-u}}{u}\,du \,,  \q z\in \CC.
     \eqno(2.6)$$
Indeed, such a function, that we refer to as the {\it modified exponential integral},
turns out to be entire, being the primitive of an entire function.
We found an account of this function in Ch. 6 of the  handbook NIST edited by Nico Temme 
\cite{NIST}, where it is referred to as the {\it complementary exponential integral}.
%The relation with the classical Exponential integrals will be given in (D.11).
%%%%%%%%%%%%%%

%%%%%%%%

%\subsection{Power series.}

The {\it power series expansion}
of $\Ein(z)$, valid in all of $\CC\,,$
can be easily obtained by term-by-term
integration and reads
$$
 \Ein(z) := z - \frac{z^2}{2 \cdot 2!} +
   \frac{z^3}{3 \cdot 3!} -  \frac{z^4}{4 \cdot 4!} +  \dots
= \sum_{n=1}^{\infty}
   (-1)^{n-1}  \frac{z^n}{ n\,n!}.\eqno(2.7)$$

%\subsection{Relation between $\,\Ein(z)\,$ and $\,\EE_1(z)= -\Ei(-z)$.}

The relation between $\,\Ein(z)\,$ and $\,\EE_1(z)=-\Ei(-z)\,$ can be obtained from
the series expansion of the incomplete gamma function $\Gamma (\alpha \,,\,z)$ in the limit as
$\alpha \to 0\,, $ as shown \eg in the appendix D of the book by Mainardi
\cite{Mainardi_BOOK2010}. We get 
$$ \EE_1(z) = -\Ei(-z) = \Gamma(0\,,z) =
   -C  - \log z + \Ein (z)\,, \;
    \eqno(2.8)
$$
with $|\hbox{arg} \,z| < \pi$,
where $C= -\Gamma^\prime (1)= 0.577215\dots$,
 denotes the Euler-Mascheroni constant,

This relation is important for understanding the analytic properties
of the classical exponential integral functions
in that it isolates the
multi-valued  part represented by the
logarithmic function from  the regular part
represented by the entire function  $\Ein(z)$ given by the
 power series in (2.7), absolutely convergent in all of $\CC\,. $

%%%%%%%%
\subsection{The Exponential integrals in the time domain}

Let us consider in the time domain $t>0$  the following two causal functions $\phi(t), \psi(t)$
related to Exponential integrals:
$$ \phi(t) := \EE_1(t)\,,  \q t>0\,, \eqno(2.9)$$
$$ \psi(t) := \Ein(t) = C +\log \,t +\EE_1(t)\,, \q t>0\,,\eqno(2.10)$$
The corresponding Laplace transforms, analytically continued in the complex $s$ plane cut along the negative real axis,   turn out to be:
$$ \L\{\phi(t)\}(s)  =  \rec{s}\, \log (s +1)\,, \q   s \in\CC^-\,,
\eqno (2.11)$$
$$ \L\{\psi(t)\} (s);  =  \rec{s}\, \log \left(\rec{s} +1\right)\,, 
\q   s \in\CC^-\,.
\eqno (2.12)$$
%%%%%%
The proof of Eq. (2.11) can be found \eg in the book by Ghizzetti and Ossicini
\cite{Ghizzetti-Ossicini_BOOK1971}, see Eq, [4.6.16], pp 104--105.
The proof of (2.12) is hereafter provided in two ways.
%%%%%%%%%
The first proof  is obtained as a  consequence
of the identity (2.10), \ie
$  \Ein(t) =   \EE_1(t) + C  + \log \,t\,, $
and the Laplace transform pair
  $$\L\{\log\, t\}(s)  = - \rec{s}\,\left[ C+ \log \,s \right]\,, \q  
 s \in\CC^-\,,  $$
whose proof is found \ie
 in  \cite{Ghizzetti-Ossicini_BOOK1971},
see Eq. [4.6.15] and  p. 104.
%%%%%%%%%%%%%%%%
The second proof is direct and instructive.
For this it is sufficient to  compute the Laplace transform
of the elementary function provided by the derivative of $\Ein(t)$, that  is,
according to a standard exercise in the {theory of Laplace
transforms},
$$ \L\left\{ \frac{ 1-\e^{-t}}{t}\right\}(s)
  = \log \left(\rec{s} +1\right)\,, \q  s \in\CC^-
   \,,$$
so that, with $f(t) = (1-\e^{-t})/t$ and
$\widetilde f(s)= \L\{f(t)\}(s)$,  
$$ \psi(t) := \Ein(t) = \int_0^t f(t')\, dt'
  \,\div \,
  \frac{\widetilde f(s)}{s}  =
  \rec{s} \log \left(\rec{s} +1\right)\,, \q   s \in\CC^-\,, $$
  in agreement with (2.12),
where we have denoted with $\div$ the   juxtaposition of a function $f(t)$
with its Laplace transform $\widetilde f(s)$.
 
%%%%%%%%%

%%%
We  outline the different asymptotic behaviours of the two functions
$\phi(t)$, $\psi(t)$ for small argument ($ t \to 0^+$) and  large argument ($t \to +\infty$)
that can be easily obtained  by the known asymptotical  expressions  for $\EE_1(z)$ and $\Ein (z)$
available \eg in \cite{NIST}. 
%% using Eqs. (2.7), (2.10) and (2.15).
However, it is instructive to derive the required asymptotic representations
by using the Karamata Tauberian theory for Laplace transforms,
see    Feller \cite{Feller_BOOK1971}, Chapter XIII.5, as pointed out in \cite{Mainardi_BOOK2010}.
We have
$$ \phi(t) \sim
\begin{cases}
\log (1/t) \,, & \q  t \to 0^+\,, \\
 {\ds \e^{-t}/t}\,, & \q  t \to +\infty\,,
\end{cases}
\eqno(2.13) $$
%%%%%%%%%
$$ \psi(t)
\sim
\begin{cases}
 t\,, & t \to 0^+\,, \\
C + \log \, t\,, &  t \to +\infty\,.
\end{cases}
\eqno(2.14) $$
We note that the modified exponential integral (2.10) was adopted by Becker for describing the creep law in his 1926 viscoelastic model
\cite{Becker_1926} and more later in 1982 by Strick and Mainardi
\cite{Strick-Mainardi_1982} and recently by Mainardi et.al. \cite{Mainardi-Masina-Spada_2019}, even if the modified expression of the exponential integral was not known to Becker himself.
The Becker model is also discussed in \cite{Gross_BOOK1953}, 
in \cite{Mainardi_BOOK2010} and references therein.
The main reason for adopting this model in linear viscoelastitcity is that the derivative of the 
 function (2.10) in the time domain $t>0$ is interpreted (unless a suitable normalization factor) as the rate of creep: 
 $$\psi'(t)
:= \frac{d}{dt} \, \Ein (t) = \frac{1-\e^{-t}}{t}=
\frac{1}{t} - \frac{\e^{-t}}{t}, \q  t>0, \eqno(2.15)$$  
It is straightforward to prove that  the function (2.15) is completely monotonic (CM) 
 (that is a non-negative function function for $ t>0$  with infinitely many derivatives with alternating sign)
 as shown hereafter. 
 Indeed  the proof is carried out by the following  steps.
 At first, based on excellent book by Schilling et al.
\cite{Schilling-et-al_BOOK2012} that is the most recent treatise on CM functions,
we recognize that $\psi'(t)$       is  CM being  a non-negative linear combination
of two CM functions, $1/t$ and $\e^{-t}/t$. 
Then, in view of the Bernstein theorem  the function  $\psi'(t)$ can be interpreted as the Laplace transform of a positive measure, that is represented by the following integrals
 %%
% \begin{equation}
%\label{CM-spectra}
$$
\frac{d\psi}{dt} (t)= \int_0^\infty \!\!
\e^{-rt}\, K(r)\, dr = 
\int_0^\infty \!\! \e^{-t/\tau}\, H(\tau)\, d\tau \,, 
%% \end{equation}
\eqno(2.16)
$$
where  $ K(r)\ge 0 $ and  $H(\tau) \ge 0 $ are  referred to as the   spectra
in frequency ($r$) and in time ($\tau=1/r$), respectively and read
 %\begin{equation}
 $$
 K(r)=
 \left\{
 \begin{array}{ll}  
1 & 0\le r <1, \\
0  & 1\le r<\infty;
\end{array}
\right.
\quad
H(\tau)=
 \left\{
 \begin{array}{ll}  
0 & 0 \le \tau <1, \\
1/\tau^2  & 1\le \tau<\infty.
\end{array}
\right.
% \end{equation}
\eqno(2.17)
$$
We note that  the frequency spectrum $K(r)$ can be determined from the Laplace transform of  
$\psi'(t)$  by the Titchmarsh formula that reads in our notation  noting   $\psi(0^+)=0$ 
%\begin{equation}
%label{K(r)}
$$
K(r) = \pm \frac{1}{\pi}\,
\left.
\Im[s\widetilde\psi (s)]
\right|_{s =r \e^{\mp i\pi}}\,,
\q \widetilde\psi (s)= \L \{\psi(t)\}(s) \q s\in \CC.
\eqno(2.18)
$$
This a consequence of the fact that the Laplace transform of $\psi'(t)$  is the iterated Laplace transform of the frequency spectrum, that is the Stieltjes transform of $K(r)$ and henceforth the Titchmarsh formula provides   the inversion of the  Stieltjes transform, see \eg the treatise by Widder
\cite{Widder_BOOK1946}.
Indeed, the Titchmarsh formula,  intended in the limit of $s$ tending  from above and below to the negative real axis $s=-r$ with $r>0$, reads from (2.12)  
$$ K(r) = 
\pm \frac{1}{\pi}\,\left. \Im[s\widetilde\psi (s)]\right|_{s =r \e^{\mp i\pi}}=
 \frac{1}{\pi} \Im\left[\log \left(-\rec{r} +1\right)\right]\,.$$
%LUCHKO has suggested to extend the TITCHMARSH formula putting 
%the explicit formula derived from taking the Imaginary part of (2.12)  
Then the function $\psi(t)$, being non negative with a CM derivative is of Bernstein type
 following the usual terminology in \cite{Schilling-et-al_BOOK2012}.
 
\section{The generalized modified Exponential integral via the Mittag-Leffler function}

\subsection{The new function $\Ein_\nu (z), \; \nu>0$}
  Let us now generalize the modified exponential integral  
  as depending on a real parameter $\nu >0 $ by considering 
    the Mittag-Leffler function
\begin{equation}
E_{\nu}(z)%%  \equiv E_{\nu, 1}(z)
 = \sum_{k= 0}^{+\infty} \frac{z^k}{\Gamma(k\nu + 1)}\,,
\quad z \in \CC,\quad \nu > 0,
\end{equation}
that is known to be an entire function for $\nu>0$ and generalize the exponential function 
$\exp(z)$ to which it reduces just  for $\nu=1$. For details on this transcendental function the reader is referred to the 2014 treatise 
by Gorenflo et al. \cite{GKMR_BOOK2014}.
This entire function can be seen as the particular case of the Mittag-Leffler in two parameters $\nu, \mu$
\begin{equation}
E_{\nu, \mu}(z)
 = \sum_{k= 0}^{+\infty} \frac{z^k}{\Gamma(k\nu + \mu)}\,,
\quad z \in \CC,\quad \nu, \mu  > 0 \,,
\end{equation}
for $\mu=1$.
Then  we propose to define for any $\nu>0$, in the cut plane $\CC^- $, 
%% cut along the negative real axis
%%%
\begin{equation}
\Ein_{\nu}(z)\! = \!\int_0^{z} \frac{1 - E_{\nu}(-u^{\nu})}{u^{\nu}}\, du
\!=\! {\ds \sum_{n= 1}^{+\infty}  \frac{(-1)^{n-1}\,
z^{\nu n-\nu +1}}{(\nu n-\nu +1)\, \Gamma(n\nu + 1)}.}
\end{equation}
For $\nu=1$   the function (2.6) is recovered, that is  we have $\Ein_1(z) = \Ein(z)$.
   
   \subsection{The limiting regularized case $\nu=0$ for $t>0$}
We note that the limiting case $\nu =0$ requires special attention because in this case the Mittag-Leffler is no longer  an entire function,
Indeed,   limiting to the Mittag-Leffler function in the non-negative time domain,
 we see from the corresponding plots  that for $\nu\to 0$ the function  becomes discontinuous around  $t=0$ assuming the value 1 at $t=0$ and 1/2 at $t>0$, so we get 
\begin{equation}
\Ein_0(t) = \frac{t}{2}, \q t>0.
\end{equation}
  However, in this limiting case, we recover {\it formally}  this result  by summing according 
  to Ces{\`a}ro the undefined series of the corresponding limit of the 
Mittag-Leffler function, known as Grandi's series
\begin{equation}
\label{Grandi's series}
E_0(-t^0) = \sum_{k=0}^\infty (-1)^k =
1-1+1-1+\cdots = \frac{1}{2}\,,\q t>0\,.
\end{equation} 
This series is a particular realization of the so called Dirichlet $\eta$
function \cite{NIST}. The latter is part of a broad class of function series, known as
Dirichlet series, more known in rheology as Prony series, that have recently found new physical applications in the so-called Bessel models, see e.g. 
\cite{Giusti-Mainardi_2016,Colombaro-Giusti-Mainardi_2017,Giusti_2017}.

\subsection{The generalized Becker model in linear viscoelasticity}
Following the analysis in the paper by Mainardi et al.
\cite{Mainardi-Masina-Spada_2019}
it is straightforward to  introduce the generalized Becker model .
The corresponding creep function and the  rate of creep
 read with  $0< \nu \le 1, \; t\ge 0$,
\begin{equation}  
 \psi_\nu(t) =
 \left\{
\begin{array}{ll} 
\Gamma(1+\nu)\,  \Ein_\nu(t) \\
 \Gamma(1+\nu)\,
 {\ds \sum_{n= 1}^{+\infty}(-1)^{n-1} 
  \frac{t^{\nu n-\nu +1}}{(\nu n-\nu +1)\, \Gamma(1+ n\nu)},}
  \end{array}
  \right.
 \end{equation}
 and 
 \begin{equation}  
 \psi'_\nu(t) := \frac{d\psi_\nu}{dt} 
= \Gamma(1+\nu)\,
 \sum_{n= 1}^{+\infty}(-1)^{n-1} 
  \frac{t^{\nu n-\nu }}{ \Gamma(1+ n\nu)},
 \end{equation}
 where the factor $\Gamma(1+\nu)$ is settled  to get $\psi'(0)=1$.
The regularized limiting result for $\nu=0$   corresponds in linear viscoelasticity to the linear creep law for a Maxwell body.
As a consequence,  our  generalized Becker model is  ranging from the Maxwell body
 at $\nu=0$ to the Becker body at 
 $\nu=1$.
 
 As in \cite{Mainardi-Masina-Spada_2019},
 in the following figure (Fig. 1) we show versus time the creep function   
 $\psi_\nu(t)$ and its derivative (the rate of creep) $\psi'(t)$ in a linear scale $0\le t \le 10$ for the particular values of $\nu=0, 0.25, 0.50, 0.75, 1$,  
from where we can note the tendency to the linear creep law for the Maxwell model as
$ \nu \to 0^+$.

\begin{figure}[h!]
\begin{center}
\includegraphics[scale=0.40]
{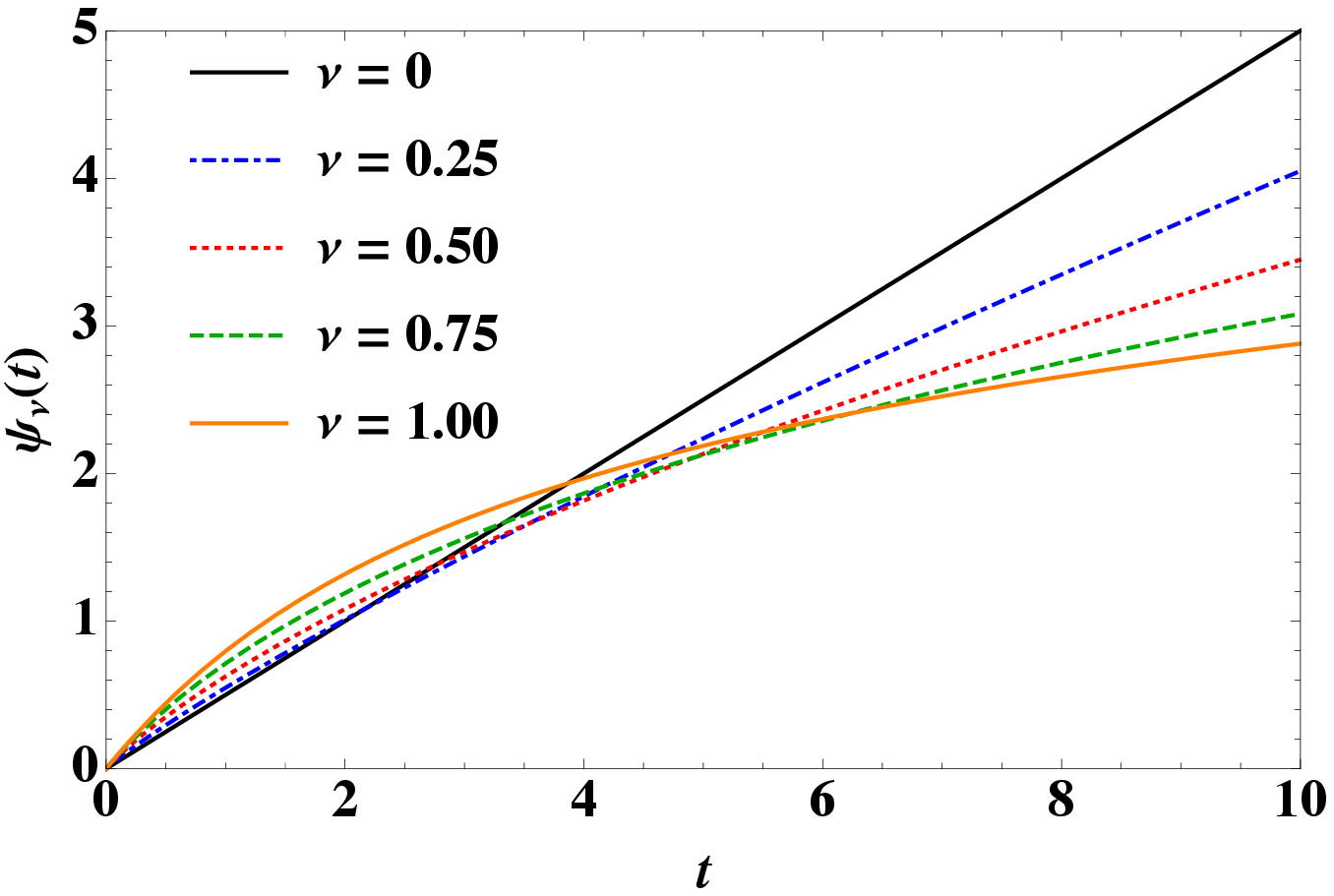}
\includegraphics[scale=0.40]
{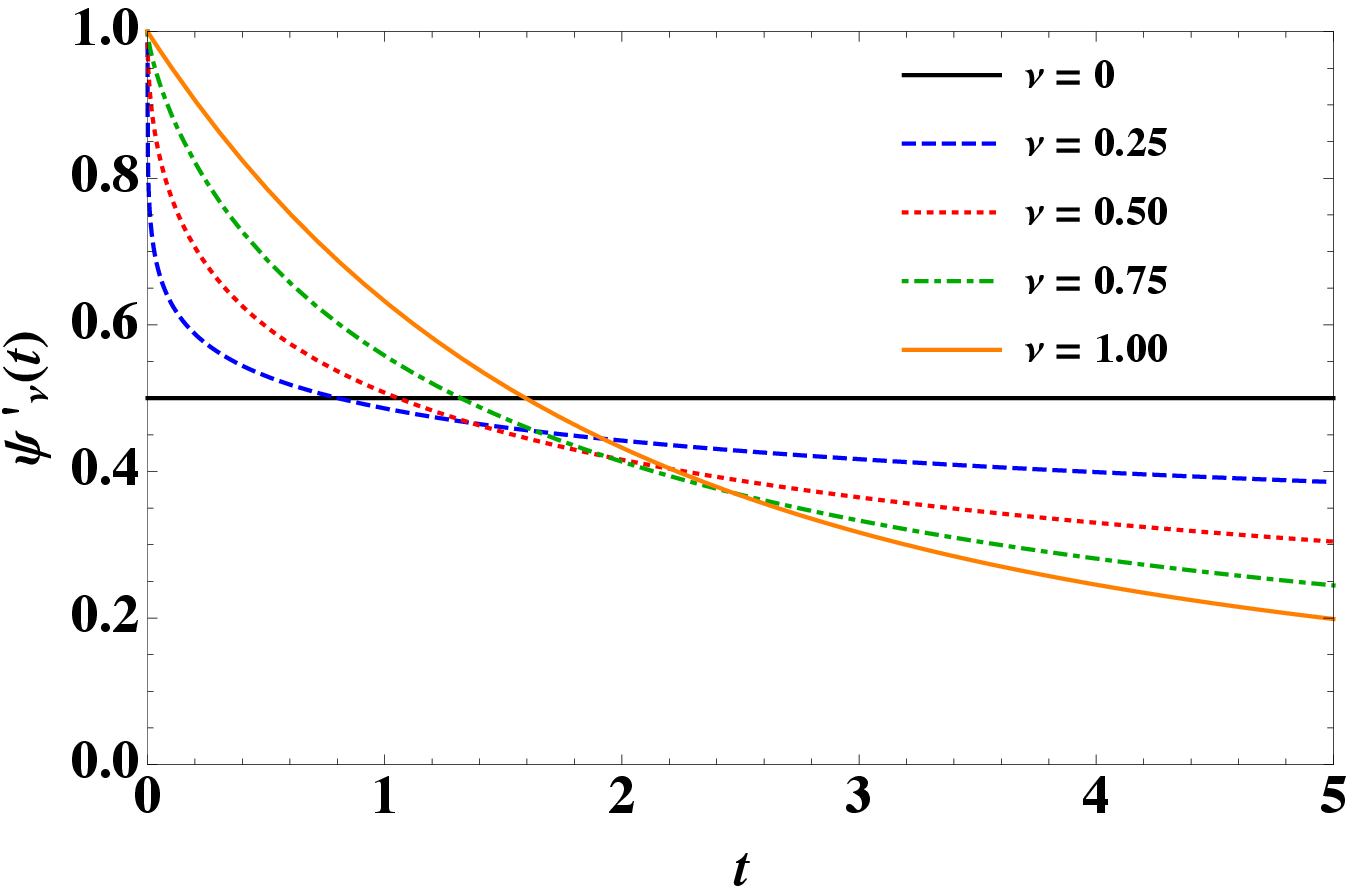}
%{CreepRate.eps}
%{PSI-LINEAR.eps}
%{Fig3.png}
 \end{center}
 \vskip-0.5truecm
\caption{The creep function $\psi_\nu(t)$ (left)
and the rate of creep  $\psi'_\nu (t)$  (right) for selected values of $\nu \in [0,1]$.}
\end{figure}

 \newpage
 
We point out  that also for $0< \nu<1$ 
the  functions characterizing the generalized Becker model 
%% the corresponding non dimensional functions
 $\psi_\nu(t)$ and 
$\psi'_\nu(t)$ keep the property to be Bernstein and CM functions as it is for the Maxwell and Becker bodies. 
The proof is straightforward because it is obtained by following the same steps of the Becker model illustrated in the previous section.This is due to the property of complete monotonicity of the 
Mittag-Leffler function  even if divided for a power law function of Bernstein type. Indeed it is known that a CM  function divided by a Bernstein function is still   CM, see \eg
\cite{Schilling-et-al_BOOK2012}.
Furthermore,  the function $\psi_\nu(t)$, being non negative with a CM derivative is of Bernstein type following the usual terminology in  \cite{Schilling-et-al_BOOK2012}.

Indeed, following the same arguments of the previous section we can compute the spectral distributions
$K_\nu(r)$, $H_\nu(\tau)$ related to $\psi'_\nu (t)$  applying the Titchmarsh formula to the Laplace transform of $\psi'_\nu(t)$ . 
 However, for $0<\nu<1$ 
 the Laplace transform of the rate of creep  is not known in analytic form  so that it can be obtained  integrating term by term the series representation of 
 the original function. 
 The series in the Laplace domain of the functions $\psi_\nu(t)$ and 
 $\psi'_\nu(t)$ turn out to be, with
 $ 0< \nu \le 1, \; s\in \CC^-$:
 \begin{equation}  
 \L\{\psi_\nu(t)\}(s) = 
 \Gamma(1+\nu)\,
 \sum_{n= 1}^{+\infty}(-1)^{n-1} 
  \frac{s^{\nu -n\nu -2}\Gamma (n\nu-\nu+2)}
  {(\nu n-\nu +1)\, \Gamma(1+ n\nu)},
  \end{equation}
 \begin{equation}  
 \L\{\psi'_\nu(t)\}(s) = 
 \Gamma(1+\nu)\,
 \sum_{n= 1}^{+\infty}(-1)^{n-1} 
  \frac{s^{\nu -\nu n-1}\Gamma(n\nu-\nu+1)}{ \Gamma(1+ n\nu)}.
 \end{equation}
 
We show in Fig. 2 the spectra in frequency and in time of the rate of the creep for selected values of $\nu$  in the range $0<\nu \le 1 $
obtained  by using the 
{$ MATHEMATICA^{\textregistered}$} tool box on the series 
of $\L\{\psi'_\nu(t)\}(s)$ entering the Titchmarsh formula (2.18) with 
$\widetilde \psi_\nu(s)$.  
%%%
 These spectra  turn out to be non-negative 
 (with  a semi-infinite support $[0, +\infty)$ except in  the Becker case 
 $ \nu=1$) that, as a matter of fact, are consistent with the CM 
property of $\psi'_\nu(t)$.
\begin{figure}[h!!]
    \centering
   \includegraphics[width=5.5cm]{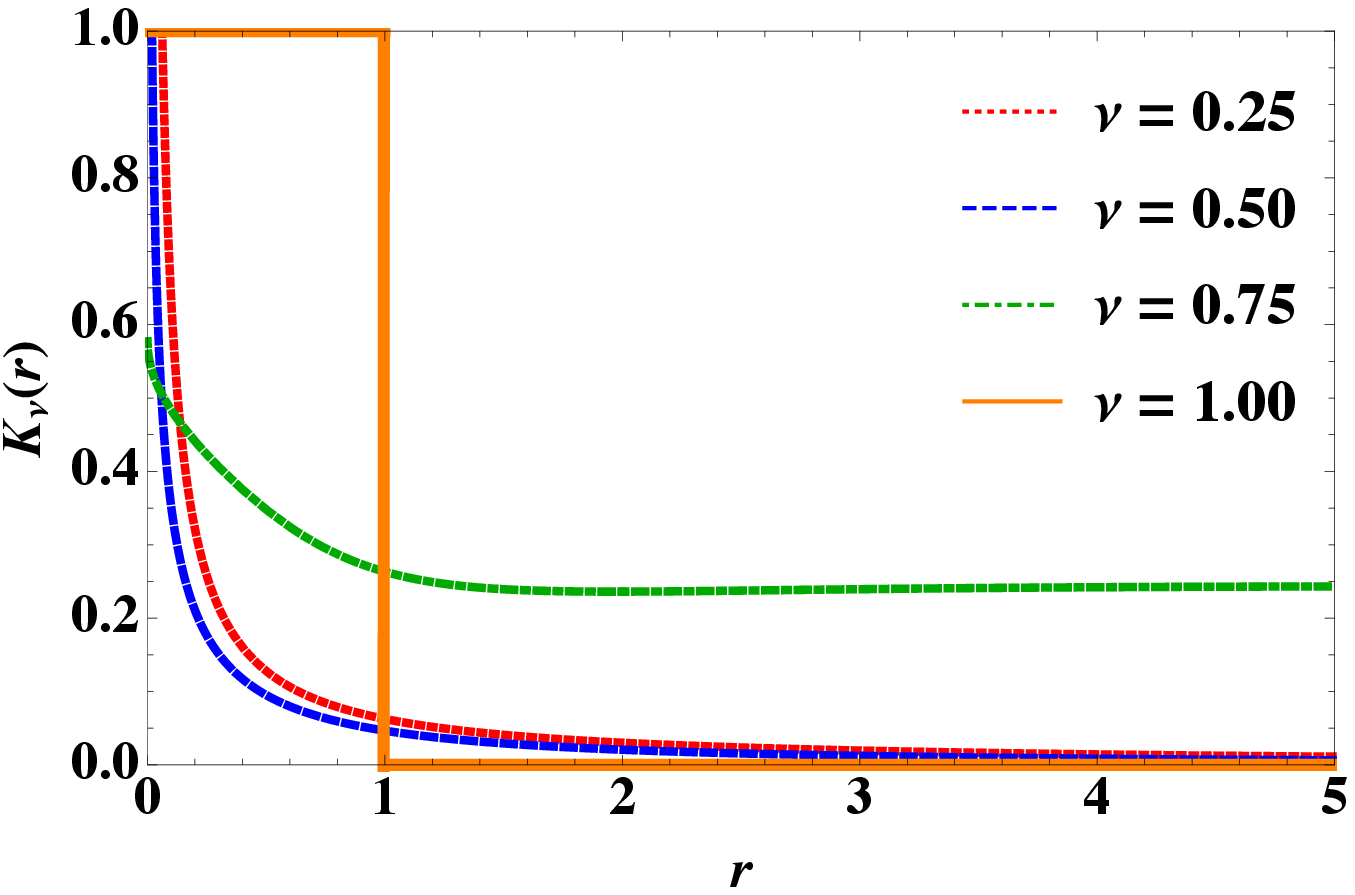}
    $\;$ 
    \includegraphics[width=5.5cm]{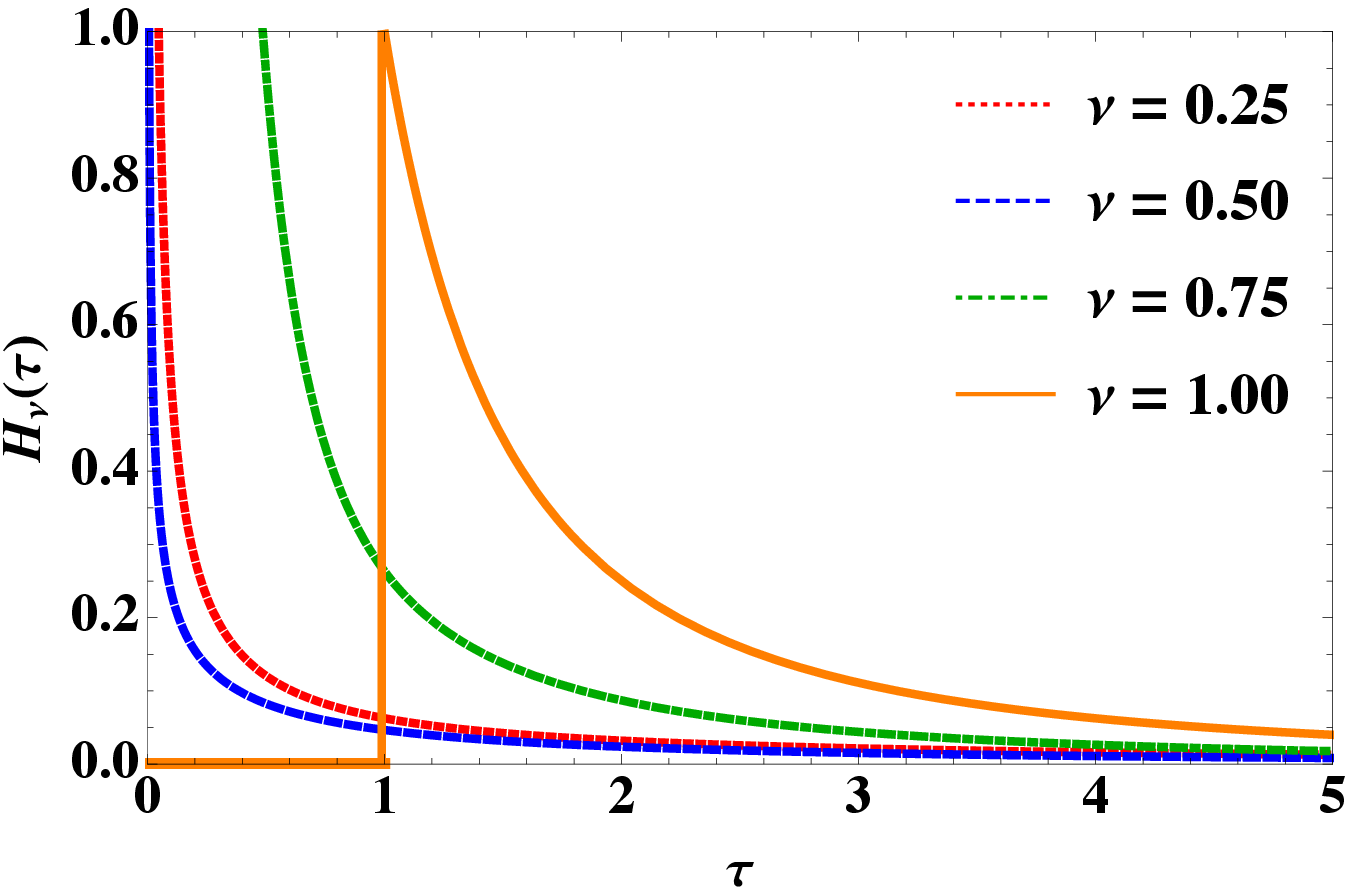}
   \vskip-0.5truecm
    \caption{The spectra for the generalized Becker model  for $\nu=0.25, 0.50, 0.75$ compared with those of the Becker model $\nu=1$: left in frequency $K_\nu(r)$; right in time $H_\nu(\tau)$.}%
\label{H-K TOTALI-spectra}%
\end{figure}
%%LUCHKO suggests to write the series expansions of the LAPLACE TRANSFORM. 

We can note
the power-law  behaviour for the spectra $K_\nu(r)$ and 
$H_\nu(\tau)$,  due to power law asymptotics of the rate of creep $\psi'_\nu(t)$  for $t\to 0$ and $t \to \infty$ 
for $0<\nu<1$, in contrast with the exponential decay of this function occurred for $\nu=1$.     
Indeed the spectrum $K_\nu(r)$ decays like  $r^{-(1+\nu)}$ at $r\to 0$ and $r \to \infty$ as we can see considering the major contribution related to the series of  $\widetilde \psi(s)$ entering the Tichmarsh formula for $s=-r$. 

\section{The  generalized Sine and Cosine integrals}
\subsection{The classical Sine and Cosine integrals}
As usual in books on special functions, see \eg \cite{NIST} we define
the {\it sine integral} and the {\it cosine integral} by the following integrals 
in the complex plane $\CC$,
\begin{equation}
\text{Si}(z) := \int_0^z \frac{\sin(t)}{t}\ \text{d} t\,, \quad
\text{Ci}(z) := -\int_{+\infty}^z \frac{\cos(t)}{t}\ \text{d} t\,,
\end{equation}
where the path does not cross the negative real axis or pass through the origin.
It is easy to recognize that 
$\text{Si}(z)$ is an (odd) entire function whereas $\text{Ci}(z)$ is a polydrome  function with a branch cut on the negative real axis.
We also recognize that following limits on the positive real axis
 \begin{equation}
 \lim_{x\to \infty} \text{Si}(x)= \frac{\pi}{2}\,, \q 
  \lim_{x\to \infty} \text{Ci}(x)= 0\,.
  \end{equation}
  In some books we find another function related to the sine integral
  \begin{equation}
  \text{si}(z) :=- \int_z^\infty \frac{\sin(t)}{t}\ \text{d} t
  = \text{Si}(z)- \frac{\pi}{2}\,.
  \end{equation}
  To remove the polydromy of $\text{Ci}(z) $
 we follow Shelkunoff and introduce the 
  (even) entire function $\text{Cin}(z)$, referred to as the {\it modified cosine integral} as
  \begin{equation}
  \text{Cin}(z) = \int_0^z \frac{1 - \cos(t)}{t}\ \text{d} t
  = \sum_{n=1}^\infty \frac{(-1)^{n-1} z^{2n}}{\Gamma(2n+1) 2n}\,, \quad z \in \CC\,.
\end{equation} 
  Because of a  known relation concerning  the exponential integrals,  we easily get the relation between the standard cosine integral and the modified cosine integral, that is
  \begin{equation}
   \text{Ci}(z) = - \text{Cin}(z) + \log z + C\,, \q z \in \CC^-\,.
   \end{equation}
   For completeness we define
\begin{equation}   
   \text{Sin}(z) =  \text{Si}(z) =
   \int_0^z \frac{\sin(t)}{t}\, \text{d} t
   = \sum_{n=1}^\infty \frac{(-1)^{n-1} z^{2n-1}}{\Gamma(2n) (2n-1)}
   \,, \quad z \in \CC\,.
   \end{equation}
   
   Now we like to generalize  the   entire functions 
   $\text{Sin}(z)$ and $\text{Cin}(z)$  in the cut-plane $\CC^-$ 
   by introducing a real (positive) parameter $\nu$ and using the Mittag-Leffler functions. We restrict our analysis to $\nu \in (0,1)$ with the request to recover the above classical functions 
 for $\nu=1$.
 We first note that in the literature there exist some different definitions of the fractional circular functions depending on the parameter $\nu \in (0,1)$ adopted  to generalize the classical circular function for $\nu=1$, see
 the 1999 textbook by Podlubny \cite{Podlubny_BOOK1999}, p.19, 
      Eqs. (1.69)-(170), (1.71)-(1.72), formerly provided by Tscytlin in 1984 \cite{Tscytlin_BOOK1984} and Luchko and Srivastava in 1995   \cite{Luchko-Srivastava_1995}.
      
\subsection{The generalized Sine and Cosine integrals}
For generalizing  Sine and Cosine integral functions in the cut plane $\CC^-$, 
we agree to start with the fractional circular functions defined more recently        
       in the 2014 textbook by Herrmann 
      \cite{Herrmann_BOOK2014}, Ch. 6.
      %%%%
   They are related   to the Mittag-Leffler functions defined in (3.2)  and read
   \begin{equation}    
\hskip-0.5truecm %      \!\! 
      \begin{cases}  %%\left\{
  \hskip-0.1truecm 
    \begin{array}{ll}
     \sin_\nu(z) = z^\nu\, E_{2\nu, 1+\nu} \left(-z^{2\nu}\right) =
  {\ds  \sum_{k = 0}^{+\infty} (-1)^k
  \, \frac{z^{(2k+1)\nu}} {\Gamma( (2k+1)\nu +1) }, 
  }        
   \\
   \cos_\nu(z) = E_{2\nu, 1} \left(-z^{2\nu}\right) =  
   {\ds  \sum_{k = 0}^{+\infty} (-1)^k
  \, \frac{z^{2k\nu}} {\Gamma( 2k\nu +1) }, 
 }       
\end{array}
  \end{cases}  %%  \right .
   \end{equation}
   For the {\it fraction sine} and {\it fractional cosine} we get  
   the following plots  in Fig. 3 in the range $0\le x \le 10$, for $0<\nu \le 1$:
       \begin{figure}[h!]
\begin{center}
\includegraphics[scale=0.40]
{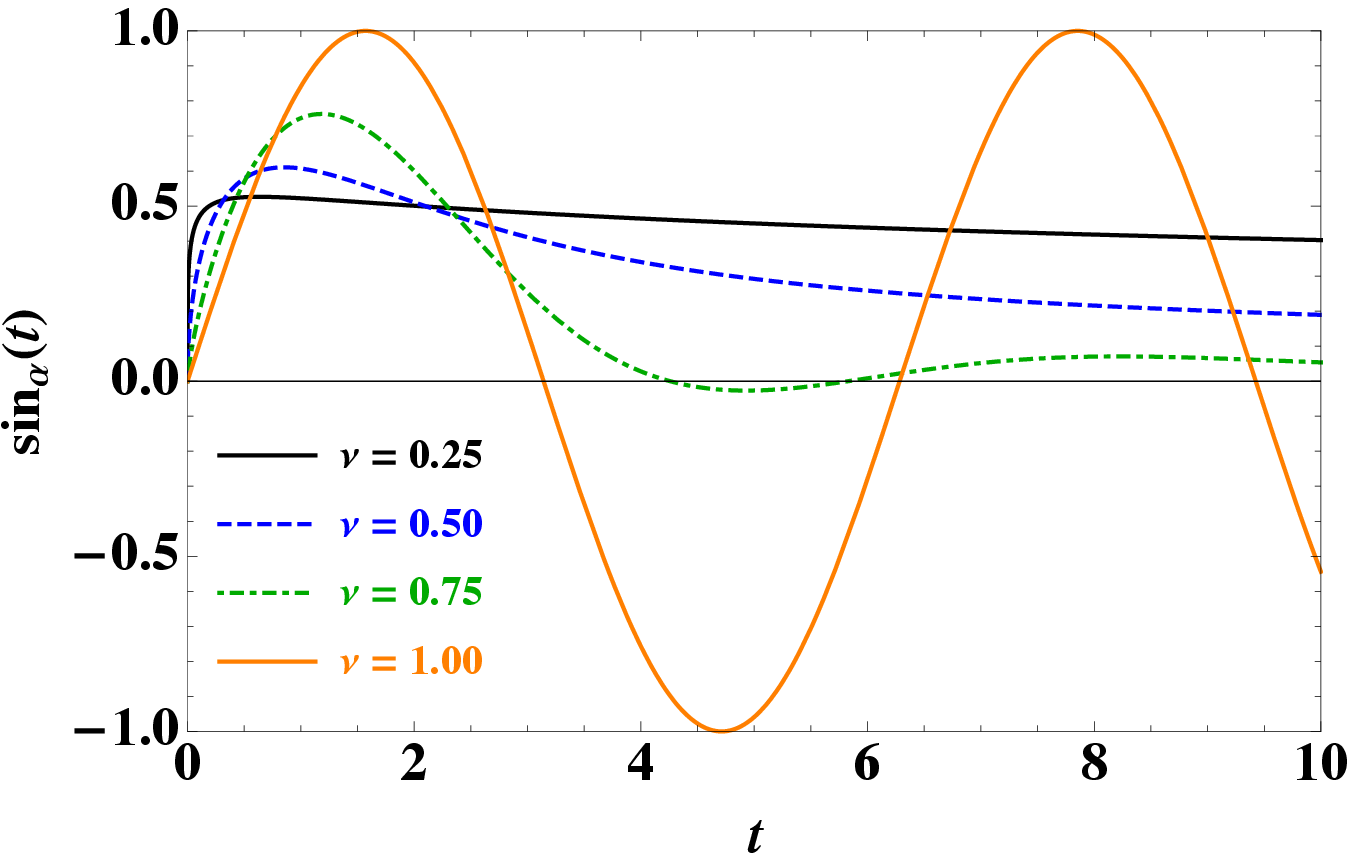}
\includegraphics[scale=0.40]
{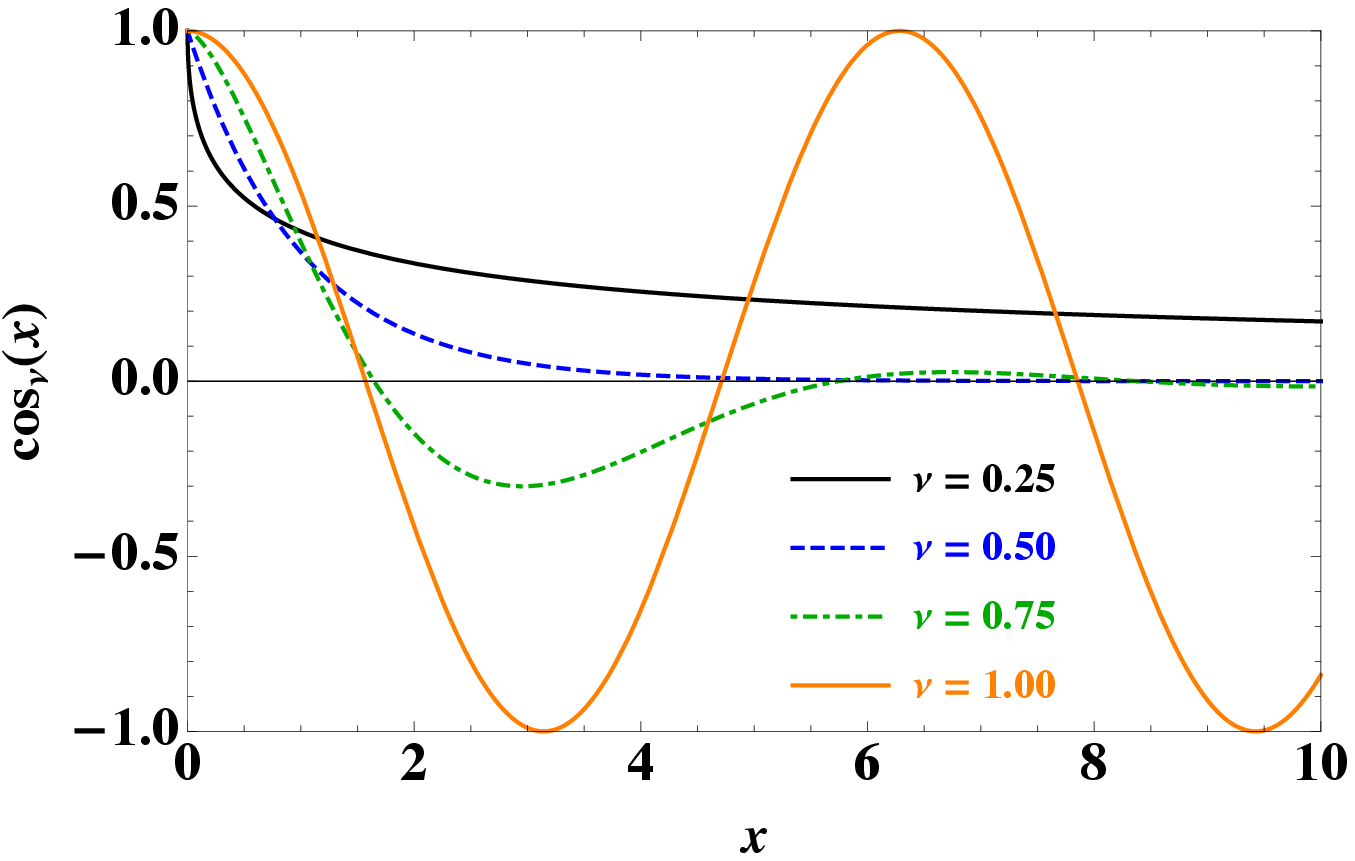}
 \end{center}
 \vskip -0.5truecm
\caption{The fractional sine  (left) and  the fractional  cosine  (right) 
  versus $x\in[0,10]$  for $\nu = 0.25, 0.50, 0.75, 1$ according to the approach by Herrmann}
\end{figure}
 We easily recognize because of the power law asymptotics of the Mittag-Leffler functions  that both the fractional circular functions
 are decaying for $\nu\in (0,1)$ as a power law after a finite number of oscillations. 
 
 Furthermore for $\nu=0$ the fractional sine and cosine  functions   are expressed by  Grandi's series so the corresponding plots tend as $\nu\to 0$ to the constant value 0.5 for $x>0$ with the value 0 and 1 at $x=0$, respectively.
   
Then we define the  functions $\text{Sin}_\nu(z)$ and $\text{Cin}_\nu (z)$
by  the following integrals 
  \begin{equation}
   \!\!\left\{
   \begin{array}{ll}
     \!\! \text{Sin}_\nu(z) =  
   {\ds \int_0^z \frac{\sin_\nu(t)}{t^\nu}  \ \text{d} t}    
   \,, \\ 
     \!\! \text{Cin}_\nu(z) = 
     {\ds \int_0^z \frac{1 - \cos_\nu(t)}{t^\nu}\ \text{d} t}       
     \,.
   \end{array}
   \right .
   \end{equation}
   
  %%%%%%%%%%%%
As a consequence we derive the following power series for the desired functions 
 \begin{equation}
   %\left\{
   \begin{cases}
\text{Sin}_\nu(z) =  
  {\ds  \sum_{k = 0}^{+\infty} (-1)^k
  \frac{z^{2k\nu+1}}{\Gamma((2k+1)\nu+1)\, (2k\nu+1) },
}   
   \\
  \text{Cin}_\nu(z) = 
     {\ds  \sum_{n = 1}^{+\infty} (-1)^{n-1} 
  \frac{z^{(2n-1)\nu+1}} {\Gamma(2n\nu+1)\, ((2n-1)\nu+1) }.
  }
   \end{cases}
   %\right .
   \end{equation}
  
      We report in Fig. 4 the plots for $z=x$ of these functions for selected values of $\nu$, that is $\nu=0.25, 0.50, 0.75, 1$, in the range $0\le x\le 10$.
      
        We note that   for $\nu=1$ the plots by Gatteschi \cite{Gatteschi_BOOK1973} of $\text{Sin}(x)$, and $\text{Cin}(x)$ are reproduced in our rage $x \in[ 0,10]$.
   \newpage
     \begin{figure}[h!]
\begin{center}
\includegraphics[scale=0.25]
{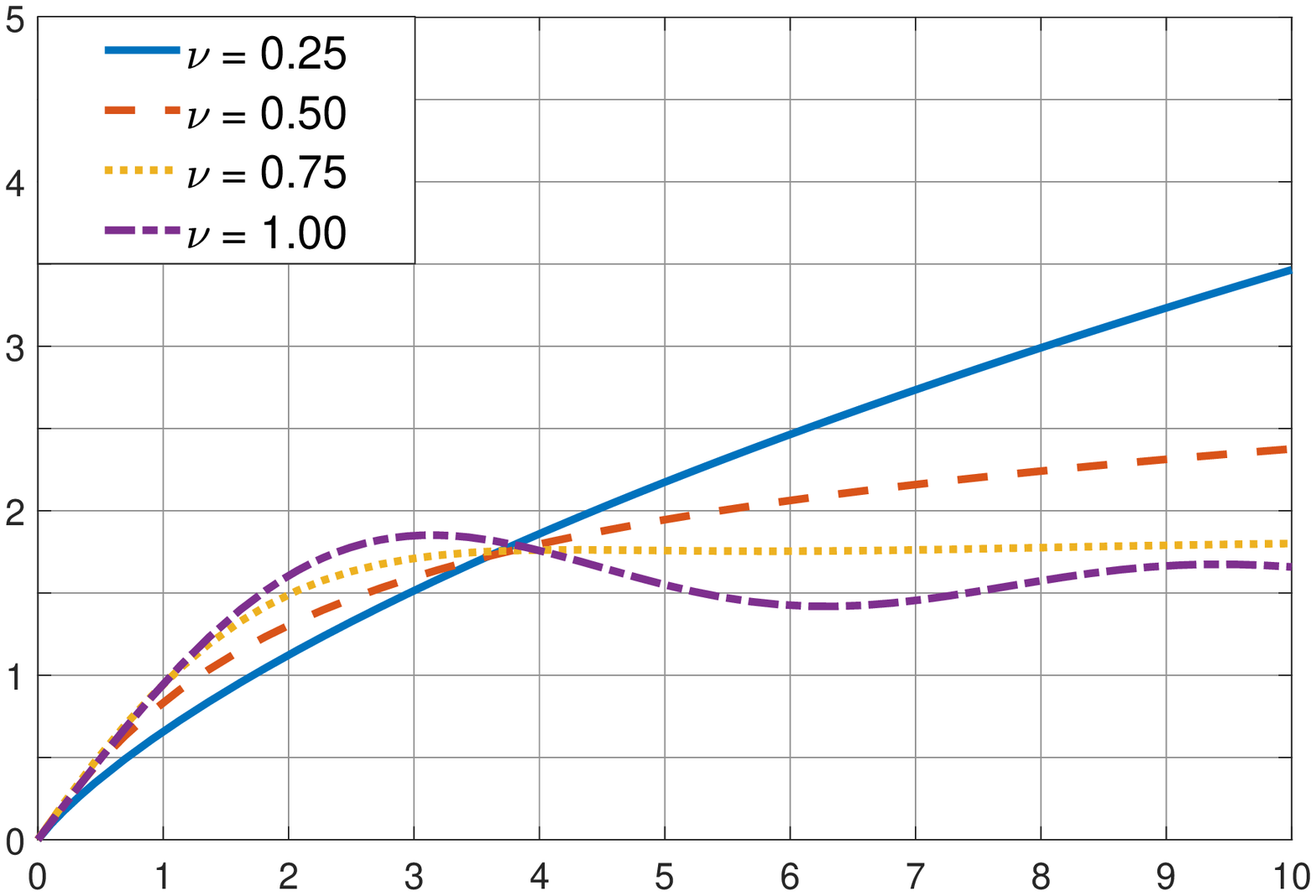}
\includegraphics[scale=0.25]
{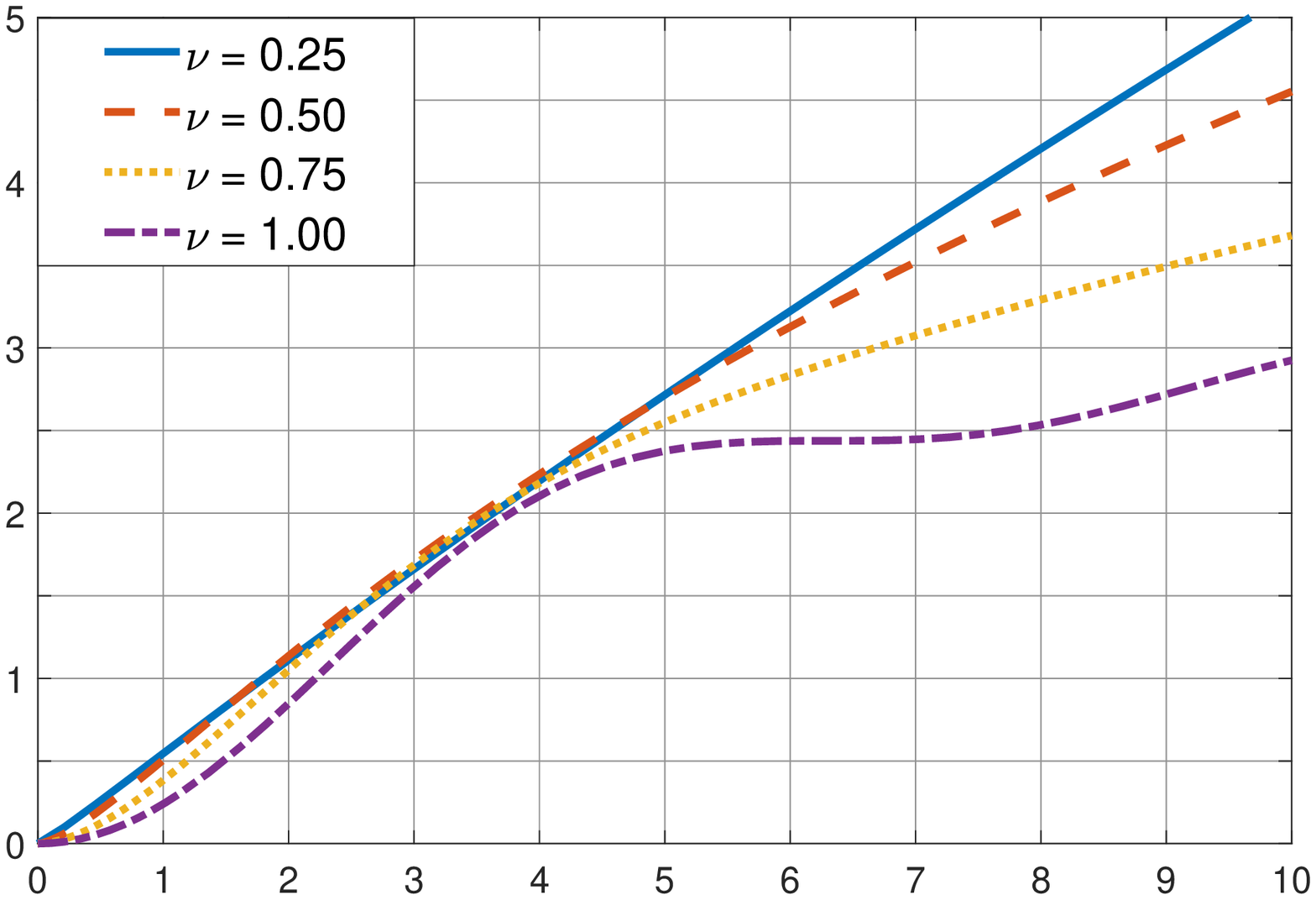}
 \end{center}
 \vskip -.05truecm
\caption{The generalized sine integral (left) and  the generalized cosine integral (right) 
  versus $x\in[0,10]$  for $\nu = 0.25, 0.50, 0.75, 1$,
based on  the approach by Herrmann for the fractional circular functions.}  
  \end{figure}
%\newpage
%%%%%%%%%%%%%%%%%%     
 \section{Conclusions}
In this paper, after having recalled  the properties of the  
Schelkunoff modification of the Exponential integral, we have  generalized it with the Mittag-Leffler function.
So doing we got a new special function (as far as we know)  that was shown to be   relevant in linear viscoelasticity because of its complete monotonicity properties in the time domain.   Indeed this new model depending on the parameter $\nu \in [0,1$] allows a  transition from 
from the standard Maxwell model for $\nu =0$  to the Becker model 
for $\nu=1$.
We have also considered  am approach  the generalized sine and cosine integral  functions
based  on a particular definition of the fractional sine and cosine functions, that appear reasonable in absence of a unique definition  for them. 
% However for these functions we have not yet found suitable applications that can  be used to justify the different definitions.
   We believe that  this point  is  still open
   %% not so relevant 
   because  also the fractional circular functions    have not found a precise application up to now, as far as we know.
        
         \section*{Acknowledgments}
	The work of the first author  (FM) has been carried out in the framework of the activities of the National Group of Mathematical Physics (GNFM, INdAM).
 The authors are  particularly grateful to Prof. Yuri Luchko and to PhD Andrea Gusti   for  helpful discussions.
 %% and to Mr Armando Consiglio for help in correcting Fig 4. 	
	
 %5 \vskip -2.5truecm

%%%%%%%%%%%%%%%%%%%%%%%%%%%%%%%%%%%%%%%%%%%%%%%%%%%%%%

%%%%End of the main text
%%%%%%%%%%%%%%%

 \end{document}